\documentclass[twoside,emlines,bezier,12pt]{article}
\usepackage{amsfonts,amssymb}

\begin{document}

\begin{center}

{\Large \bf The only regular inclines are distributive lattices}

$ $

\bigskip {Song-Chol Han$\textrm{*}$,\quad  Hak-Rim Ri}

\vskip 3mm

{\small \centerline{{\it Faculty of Mathematics, Kim Il Sung
University, Pyongyang, DPR Korea}}}

\renewcommand{\thefootnote}{\alph{footnote}}
\setcounter{footnote}{-1} \footnote{$\textrm{*}$ Corresponding
author.}
\end{center}

\vspace{0.3cm}

{\small \noindent {\bf Abstract}
\medskip

An incline is an additively idempotent semiring in which the
product of two elements is always less than or equal to either
factor. This paper proves that the only regular inclines are
distributive lattices, which also implies that there is no
noncommutative regular incline.
\medskip

\noindent {\it AMS subject classification:}\quad  16Y60
\medskip

\noindent {\it Keywords:}\quad  Semiring; regularity; idempotence;
incline; distributive lattice}
\medskip

$ $

{\section*{\noindent \bf 1. Introduction}}

Inclines are additively idempotent semirings in which products are
less than or equal to either factor. The study of inclines is
generally acknowledged to have started by Z.Q. Cao in a series of
his papers in the first half of 1980's (see \cite{Kim2004}).
Nowadays, one may clearly notice a growing interest in developing
the algebraic theory of inclines and their numerous significant
applications in diverse branches of mathematics and computer
science.

Recently, Meenakshi et al. \cite{Meen2010r} proved that an incline
is regular if and only if it is multiplicatively idempotent and
that every commutative regular incline is a distributive lattice.
Furthermore, Meenakshi et al. \cite{Meen2010g}-\cite{Meen} studied
generalized inverses, subtractive ideals, homomorphism theorems
and quotient inclines in the setting of regular inclines. But one
can easily find out that any example of noncommutative regular
inclines has never been given in \cite{Meen2010r}-\cite{Meen}.

The purpose of this paper is to check if there exist any
noncommutative regular inclines. This paper proves that every
regular incline is commutative, hence the only regular inclines
are distributive lattices.
\medskip

$ $

{\section*{\noindent \bf 2. Preliminaries}}

We recall some known definitions and facts.

A semiring is a nonempty set $R$ together with two operations of
addition $+$ and multiplication $\cdot$ such that $(R,+)$ is a
commutative semigroup, $(R,\cdot)$ is a semigroup and
multiplication distributes over addition from either side. A
semiring $R$ is said to be commutative if $ab=ba$ for all $a,b\in
R$ (see \cite{Golan1999}).

For a semiring $R$, an element $a\in R$ is said to be additively
(resp. multiplicatively) idempotent if $a+a=a$ (resp. $a^2=a$).
$R$ is said to be additively (resp. multiplicatively) idempotent
if every element in $R$ is additively (resp. multiplicatively)
idempotent.

For a semiring $R$, an element $a\in R$ is said to be regular if
$axa=a$ for some $x\in R$. $R$ is said to be regular if every
element in $R$ is regular.

An incline is an additively idempotent semiring $R$ satisfying
$x+xy=x+yx=x$ for all $x,y\in R$. On an incline $R$, a partial
order relation $\leqslant$ is defined by $x\leqslant
y\Leftrightarrow x+y=y$ for $x,y\in R$. Then $xy\leqslant x$ and
$yx\leqslant x$ hold for all $x,y\in R$. For any $x,y,z\in R$,
$y\leqslant z$ implies $xy\leqslant xz$ and $yx\leqslant zx$ (see
\cite{Cao1984}). For an incline $R$, $I(R)$ denotes the set of all
multiplicatively idempotent elements in $R$, i.e., $I(R)=\{a\in
R\mid a^2=a\}$.
\medskip

\noindent {\bf Lemma 2.1 \cite{Meen2010r}.}\quad  Let $R$ be an
incline and $x$ an element in $R$.
\medskip

\noindent (1)\quad  $x$ is regular if and only if it is
multiplicatively idempotent.
\medskip

\noindent (2)\quad  $R$ is regular if and only if it is
multiplicatively idempotent.
\medskip

$ $

{\section*{\noindent \bf 3. Main results}}

We shall use Lemma 2.1 almost everywhere in this section.
\medskip

\noindent {\bf Lemma 3.1.}\quad  Let $R$ be an incline and $x$ a
regular element in $R$. If $y\in R$, then the following hold.
\medskip

\noindent (1)\quad  $xy$ is regular if and only if $yxy=xy$.

\noindent (2)\quad  $yx$ is regular if and only if $yxy=yx$.

\noindent (3)\quad  Both $xy$ and $yx$ are regular if and only if
$xy=yxy=yx$.
\medskip

\noindent {\it Proof.}\quad  Note that $y=y+xyx$. In fact,
$xyx=(xy)x\leqslant xy\leqslant y$.

\noindent (1)\quad  If $xy$ is regular, then
$yxy=(y+xyx)xy=yxy+xyx^2y=yxy+xyxy=yxy+(xy)^2=yxy+xy=xy$.
Conversely, if $yxy=xy$, then $(xy)^2=x(yxy)=x(xy)=x^2y=xy$, thus
$xy$ is regular.

\noindent (2)\quad  If $yx$ is regular, then
$yxy=yx(y+xyx)=yxy+yx^2yx=yxy+yxyx=yxy+(yx)^2=yxy+yx=yx$.
Conversely, if $yxy=yx$, then $(yx)^2=(yxy)x=(yx)x=yx^2=yx$, thus
$yx$ is regular.

\noindent (3)\quad  follows from parts (1) and (2). \hfill
$\square$
\medskip

\noindent {\bf Lemma 3.2.}\quad  Let $R$ be an incline and $x,y$
two regular elements in $R$. Then both $xy$ and $yx$ are regular
if and only if $xy=yx$.
\medskip

\noindent {\it Proof.}\quad  ({\it Necessity})\quad  follows from
Lemma 3.1(3).

\noindent ({\it Sufficiency})\quad  If $xy=yx$, then
$(xy)^2=x(yx)y=x(xy)y=x^2y^2=xy$, thus $xy$ and $yx$ are regular.
\hfill $\square$
\medskip

\noindent {\bf Corollary 3.3.}\quad  Every regular incline is
commutative.
\medskip

\noindent {\bf Lemma 3.4.}\quad  If $R$ is a commutative incline
and $I(R)\neq \emptyset$, then $I(R)$ is a distributive lattice
with respect to the incline operations.
\medskip

\noindent {\it Proof.}\quad  For any $x,y\in I(R)$, we have
$(x+y)^2=x^2+yx+xy+y^2=x+yx+xy+y=x+y$ and $(xy)^2=xyxy=x^2y^2=xy$.
Hence $x+y\in I(R)$, $xy\in I(R)$ and $x+y=\sup\{x,y\}$ in $I(R)$.
If $x\geqslant z$ and $y\geqslant z$ for some $z\in I(R)$, then
$xy\geqslant z^2=z$, which shows that $xy=\inf\{x,y\}$ in $I(R)$.
Thus $I(R)$ is a distributive lattice. \hfill $\square$
\medskip

Lemma 3.4 is a slight modification of Lemma 2.3 in \cite{Han2004}.
\medskip

\noindent {\bf Theorem 3.5.}\quad  An incline $R$ is regular if
and only if it is a distributive lattice with respect to the
incline operations.
\medskip

\noindent {\it Proof.}\quad  If $R$ is regular, then it is
commutative by Corollary 3.3. By Lemma 3.4, $I(R)$ is a
distributive lattice with respect to the incline operations,
provided $I(R)\neq\emptyset$. By Lemma 2.1(2), $R$ is
multiplicatively idempotent and so $R=I(R)$. Thus $R$ is a
distributive lattice. The converse is obvious. \hfill $\square$
\medskip

The following example shows that in a noncommutative incline, the
product of two regular elements is not necessarily regular and its
regularity also depends on the arrangement of factors.
\medskip

\noindent {\bf Example 3.6.}\quad  Let $R=\{a,b,c,d,f\}$ be a set
of five distinct elements and define two operations $+$ and
$\cdot$ on $R$ as follows.
\medskip
\[
\begin{array}{cc}
\begin{tabular}{c|ccccc}
+& $a$& $b$& $c$& $d$& $f$\\
\hline
$a$& $a$& $a$& $a$& $a$& $a$\\
$b$& $a$& $b$& $a$& $b$& $b$\\
$c$& $a$& $a$& $c$& $c$& $c$\\
$d$& $a$& $b$& $c$& $d$& $d$\\
$f$& $a$& $b$& $c$& $d$& $f$\\
\end{tabular}
\qquad &
\begin{tabular}{c|ccccc}
$\cdot$& $a$& $b$& $c$& $d$& $f$\\
\hline
$a$& $a$& $b$& $c$& $d$& $f$\\
$b$& $b$& $b$& $d$& $d$& $f$\\
$c$& $c$& $f$& $c$& $f$& $f$\\
$d$& $d$& $f$& $d$& $f$& $f$\\
$f$& $f$& $f$& $f$& $f$& $f$\\
\end{tabular}
\end{array}
\]
\medskip

By computation, one can easily see that $R$ is a noncommutative
incline, in which both $b$ and $c$ are regular since they are
multiplicatively idempotent. But the product $bc=d$ is not regular
because $d^2=f\neq d$, while $cb=f$ is regular because $f^2=f$.
\medskip

\vspace{0.8cm}

{\renewcommand\baselinestretch{1}
\renewcommand\refname}

\end{document}